\magnification=1200
\font\title=cmssdc10 at 15pt

\def\la{\lambda}
\centerline{\title Will the announced influenza pandemic really happen?}
\bigskip
\footnote{}{Key words and phrases: influenza, pandemic, probability, model}

\bigskip
Rinaldo B. Schinazi

University of Colorado at Colorado Springs

email: rschinaz@uccs.edu
\bigskip
{\bf Abstract. }{\sl We propose two simple probability models to compute the probability of an influenza pandemic. Under a random walk model the probability that all pandemics between times 0 and 300 occur by time 150 is 1/2. Under a Poisson model with mean interarrival time of 30 years the probability that no pandemic occurs during at least 60 years is 14\%. These probabilities are much higher than generally perceived. So yes the next influenza pandemic will happen but maybe much later than generally thought.
}

\bigskip
{\bf 1. Introduction} The recent outbreaks of avian influenza have greatly increased the general anxiety about a possible new influenza pandemic. 
This is so because it is generally believed that pandemic viruses arise as reassortants of human and avian virus strains, see Holmes et al. (2005) and Potter (2001). 
Even more worrisome, the worst influenza pandemic on record (the 1918 one) may have been an entirely avian-like virus, see Tautenberger et al. (2005). While there have been several studies on how the next pandemic will arise and on what would be the best strategy to fight it (see Mills et al. (2006), Balicer at al. (2005), Meltzer et al. (1999), Osterholm (2005)), we have found no attempt to compute the probability of the next pandemic. This is not surprising given that pandemics are rather rare and we have few data points to work with. In this note we propose two very simple stochastic processes to model the occurences of pandemics. According to our models, the probability of long periods without pandemics is rather larger than generally thought. 
\bigskip
{\bf 2. The model}
We now introduce the model. Let $(X_i)_{i\geq 1}$ be a sequence of positive random variables that are independent and identically distributed. Let $S_0=0$ and
$$S_n=\sum_{i=1}^n X_i\hbox{ for }n\geq 1.$$
Finally, let $N(t)$ be defined by
$$N(t)=\max\{n\geq 1: S_n\leq t\}.$$
In words, the $X_i$ models the interarrival time between the $(i-1)$-th and $i$-th pandemic. The random variable $S_n$ is the time of the $n$-th pandemic and $N(t)$ is the number of pandemics up to time $t$. The pandemics since 1700 are well documented so we take the year 1700 to be our origin in time.
 
The time between two pandemics varies from 10 years (1889-1900 and 1958-1969) to some 50 years (1729-1781). Hence, the assumption of independent identically distributed interarrival times seems plausible.    

We will now consider two specific distributions for interarrival times.
\bigbreak
{\bf 2.1 Poisson processes.} We take as origin of time the year 1700. 
For a Poisson process the interarrival times $X_i$ are independent and exponentially distributed with rate $\lambda$. That is,
$$P(X>x)=\exp(-\la x)\hbox{ for all }x\geq 0$$
and the mean time between two pandemics is $1/\la$. 
Hence, the probability that the time between two pandemics be at least twice the mean time $1/\la$ is
$$P(X>2/\la)=\exp(-2)\sim 14\%.$$
The probability that the time between two pandemics be at least three times the mean is $\exp(-3)\sim 5\%$. 

Since 1700 there have been about 10 pandemics (the definition of pandemic is not clear cut and experts seem to disagree on the status of a couple occurences, see Potter (2001)). Hence, the average interarrival time for this small sample ($n=10$) is 30 years (an approximate 90\% confidence interval for the mean interarrival time is (14,46)). If we use $\la=1/30$ the number of pandemics between times 0 and 300 is a poisson random variable with mean 10. Figure 1 below shows the distribution of such a random variable.
\medskip
\input epsf
\epsfbox {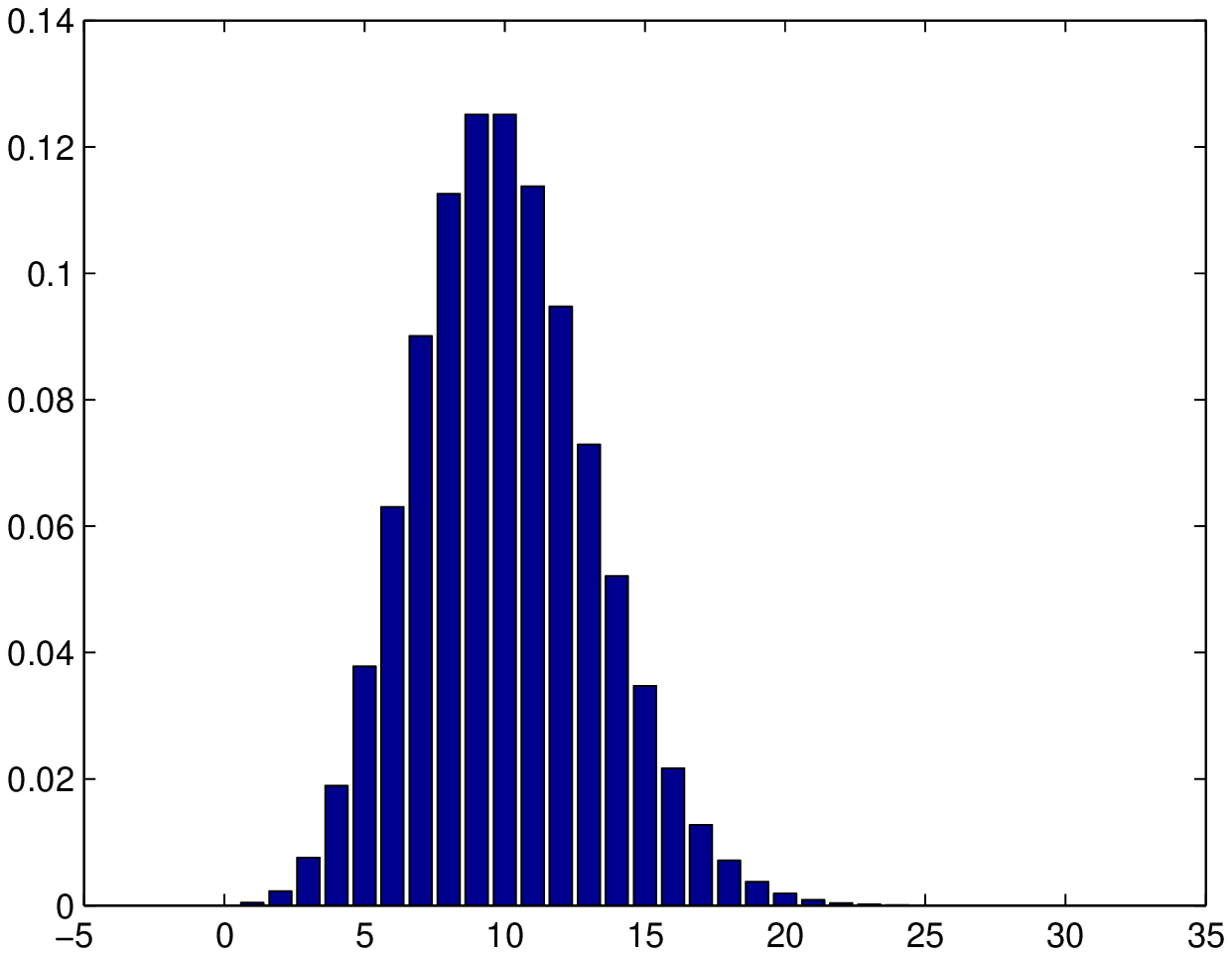}
Figure 1. This is the distribution of the number of pandemics during a 300 years period under the Poisson model with mean interarrival time of 30 years. 
\bigbreak
{\bf 2.2 Random walk.} Consider a simple symmetric random walk on the integers. If the walker is at site $i$ (where $i$ can be a positive or negative integer) at time $n$ (where $n$ is a positve integer) its position is represented by the coordinates $(i,n)$. Flip a fair coin, if it shows tails the walker jumps to $(i+1,n+1)$, if it shows heads the walker jumps to $(i-1,n+1)$. Start the walk at $(0,0)$. Let the random variable $X_j$ be the time of the $j$th return to 0 for the first coordinate. The random variables $X_j$ are independent and identically distributed. This defines a discrete time renewal process and we think of each return to 0 as being the occurence of a pandemic. In contrast to the Poisson process above we have for this model that the mean of each $X_j$ is infinite. This makes the prediction of any given $X$ very uncertain.  More precisely, we have that
$$P(X>2n)\sim{1\over\sqrt{n\pi}},$$
see (2.4) and (3.1) in Chapter  III in Feller (1968). Hence, the probability that the time between two pandemics is at least 50 years is 11\% but the probability that this time is at least 100 years is 8\%, not much less!  Observe also that the mean number of pandemics under this model is about 13 which is larger but still consistent with the 10 observed pandemics since 1700. 
\input epsf
\epsfbox {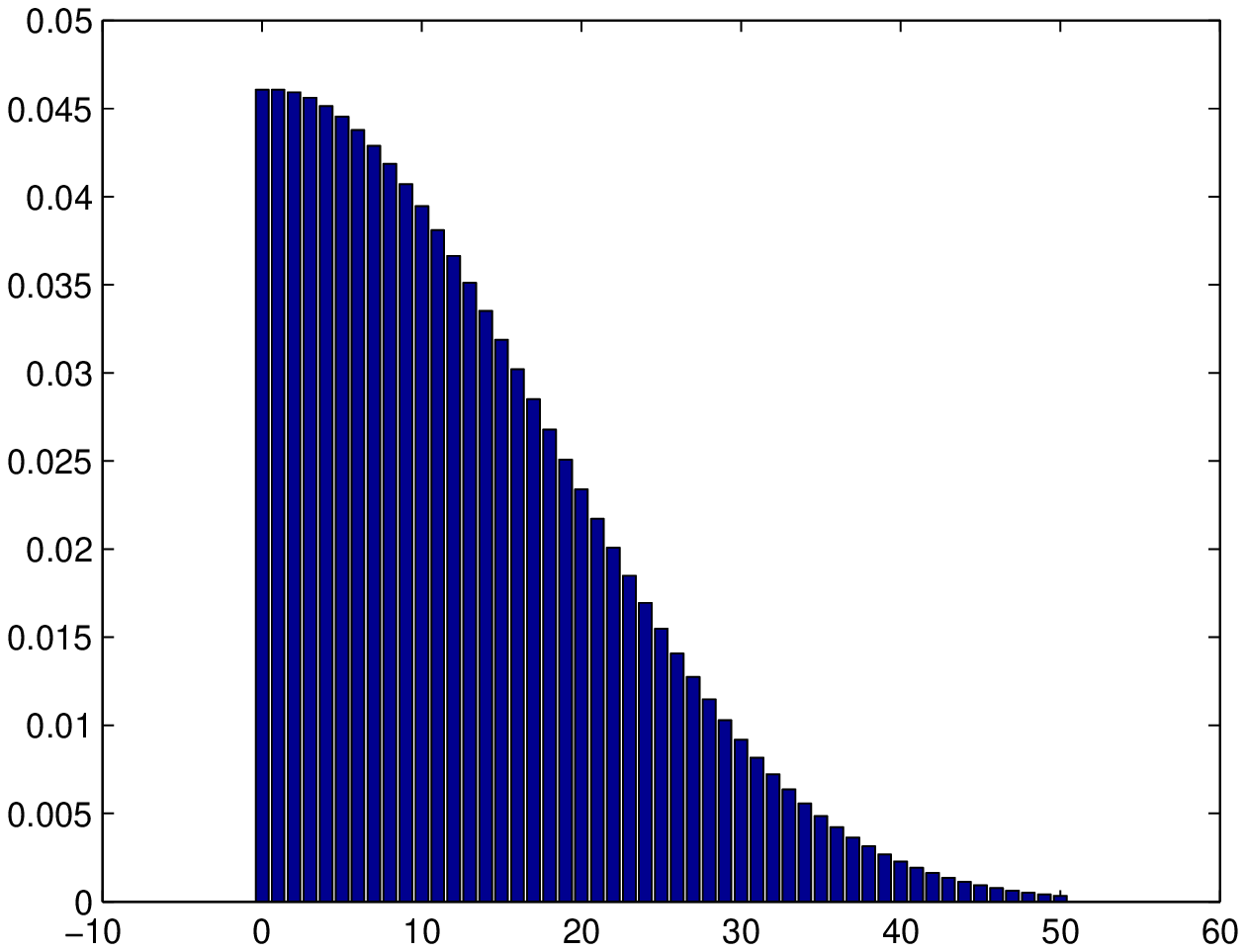}
Figure 2. This is the distribution of the number of pandemics during a 300 units period under the random walk model. This distribution can be computed exactly, see Problem 10 in III. 10 in Feller (1968). Observe how spread out this distribution is.
\input epsf
\epsfbox {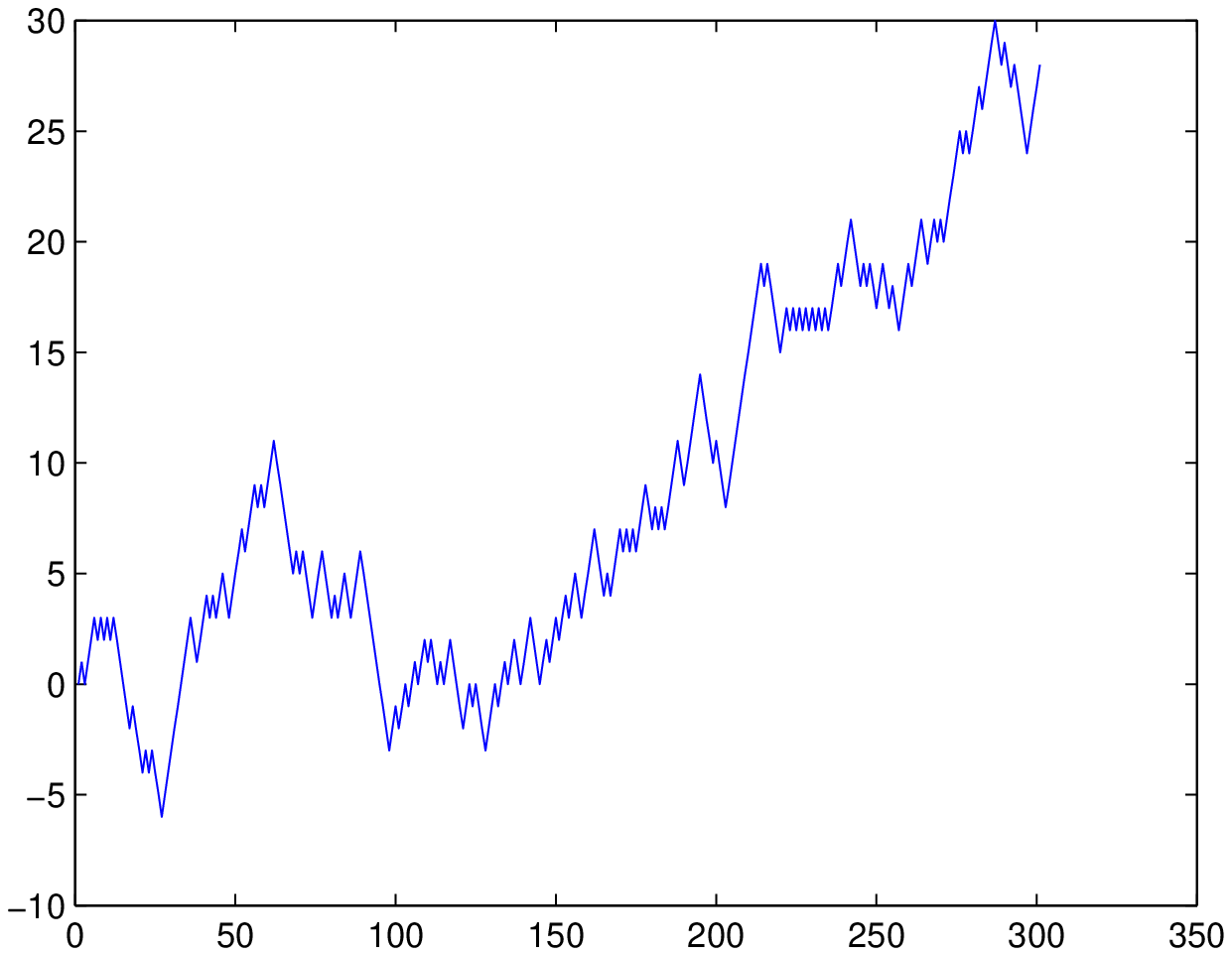}
Figure 3. This illustrates a 'typical' path for a symmetric random walk over 300 time units. In this case, there are 16 returns to 0 (corresponding to 16 pandemics). There are no pandemics between times 144 and 300. In fact, the probability that all pandemics between times 0 and 300 occur by time 150 is 1/2, see III.4 of Feller (1968).
\bigbreak
{\bf 3. Discussion} As far as we can tell this is a first attempt to compute the probability of an influenza pandemic. The occurence of a pandemic seems to be hopelessly complex in the sense that it depends on a multitude of factors, some known others unknown. Factors mentioned in the literature go from the stability of the current influenza strains to the number of pigs in China! Hence, a model taking into account all known factors would probably be as complex as the phenomenon itself, would not be that accurate and would therefore be useless. Our approach is drastically different. We do not attempt to incorporate into the model any of the factors that are believed to provoke a pandemic. Instead we treat the occurence of a pandemic as a purely random phenomenon. Among the many possible stochastic processes we picked two that are particularly random: a Poisson process and a random walk. We believe they are good starting points because the assumptions under which they are based are clearly stated and can be compared to the data and discussed. Both models are so called renewal processes. The main difference is that the time between two pandemics has a finite mean for a Poisson process but has an infinite mean for the random walk model. The latter model is admittedly quite arbitrary but it illustrates the fact that rather natural phenomena may give rise to infinite means and hence wild fluctuations.  With so few sample points it is of course difficult to assess the goodness of fit of our models but at the very least they seem plausible. 

Our main result is that, according to our models, there may be very long periods without pandemics. Under the random walk model the probability, in any $2n$ years period, that all pandemics will happen in the first $n$ years is 50\% (for any $n$). Under the Poisson model with mean interarrival time of 30 years the probability that no pandemic occurs during at least 60 years is 14\%. These probabilities seem much higher than generally perceived. Our main conclusion is that the next pandemic might not be as imminent as generally thought.
\bigbreak
{\bf References.}

R.D.Balicer, M.Huerta, N. Davidovitch and I. Grotto (2005). Cost-benefit of stockpiling drugs for influenza pandemic. Emerging infectious diseases, 11, 1280-1282.

W. Feller (1968). {\it An introduction to Probability Theory and its applications,} Third Edition, John Wiley and Sons.

E.C. Holmes, J.K. Tautenberger and B.T.Grenfell (2005). Heading off an influenza pandemic. Science, 309, 989. 

M.I.Meltzer, N.J.Cox and K.Fukuda (1999). The economic impact of pandemic influenza in the United States: priorities for intervention. Emerging infectious diseases, 5, 559-671.

C. E. Mills, J.M. Robins, C.T. Bergstrom and M. Lipsitch (2006). Pandemic influenza: risk of multiple introductions and the need to prepare for them. PLoS Medicine, 3, www.plosmedicine.org

M.T. Osterholm (2005). Preparing for the next pandemic. Foreign Affairs, July-August 2005, http://www.foreignaffairs.org

C.W. Potter (2001). A history of influenza. Journal of applied Microbiology, 91, 572-579.

R.B.Schinazi (1999). {\it Classical and spatial stochastic processes.} Birkhauser.

J.K. Tautenberger, A.H. Reid, R.M. Lourens, R.Wang, G. Jin and T.G.Fanning (2005). Characterization of the 1918 influenza virus polymerase genes. Nature, 473, 889-893.
\bye